\documentclass[12pt]{amsart}
\usepackage{amsmath}
\usepackage{amssymb}
\usepackage{amsthm}
\usepackage{amscd}

\usepackage{color}
\usepackage{hyperref}

\usepackage{tikz-cd}

\usepackage[latin1]{inputenc}
\usepackage{mathpazo}
\usepackage[scaled=.95]{helvet}
\usepackage{courier}

\swapnumbers
\headsep=1cm \numberwithin{equation}{section}
\usepackage[colorinlistoftodos]{todonotes}

\newtheorem{theorem}{Theorem}[section]
\newtheorem{lemma}[theorem]{Lemma}
\newtheorem{proposition}[theorem]{Proposition}
\newtheorem{corollary}[theorem]{Corollary}

\theoremstyle{definition}
\newtheorem{definition}[theorem]{Definition}
\newtheorem{remark}[theorem]{Remark}
\newtheorem{remark and definition}[theorem]{Remark and Definition}
\newtheorem{remark and notation}[theorem]{Remark and Notation}
\newtheorem{convention}[theorem]{Convention}

\newtheorem{notation}[theorem]{Notation}
\newtheorem{example}[theorem]{Example}

\newtheorem{conjecture}[theorem]{Conjecture}
\newtheorem{question}[theorem]{Question}
\newtheorem{questions}[theorem]{Questions}

\newcommand\Hom{\operatorname{Hom}}
\newcommand\Ext{\operatorname{Ext}}
\newcommand\Tor{\operatorname{Tor}}

\newcommand\depth{\operatorname{depth}}

\newcommand\Ker{\operatorname{\Ker}}
\newcommand\Coker{\operatorname{Coker}}

\newcommand\pd{\operatorname{pd}}
\newcommand\id{\operatorname{id}}

\newcommand\Der{\operatorname{Der}}

\newcommand{\m}{\mathfrak{m}}

\begin{document}
\dedicatory{Dedicated to the memory of  J\"urgen Herzog }

\title[Some Homological Conjectures Over Idealization Rings]{Some Homological Conjectures Over Idealization Rings}

\date{}

\author{I. J. Nascimento}
\address{Universidade de S{\~a}o Paulo -
ICMC, Caixa Postal 668, 13560-970, S{\~a}o Carlos-SP, Brazil}
\email{igorjdnascimento@gmail.com}

\author{V. H. Jorge-P\'erez}
\address{Universidade de S{\~a}o Paulo -
ICMC, Caixa Postal 668, 13560-970, S{\~a}o Carlos-SP, Brazil}
\email{vhjperez@icmc.usp.br}

\author{T. H. Freitas}
\address{Universidade Tecnologica Federal do Paran\'a, 85053-525, Guarapuava-PR, Brazil}
\email{freitas.thf@gmail.com}

\date{\today}
\thanks{}

\keywords{}
\subjclass[2010]{13D22, 13D07, 13H05, 13H10, 13C10}

\thanks{Corresponding author: I. J. Nascimento}

\begin{abstract} Let $(R,\m,k)$ be a Noetherian local ring and let $M$ be a finitely generated $R$-module. The main focus of this paper is to give positive answers for  some long-standing homological conjectures over the idealization ring $R\ltimes M$. First, if $N$ is a $R\ltimes k$-module, we show that the vanishing of $\Ext_{R\ltimes k}^{i}(N,N\oplus (R\ltimes k))$ for some $i\geq 3$ gives that  $N$ is free, and this provides  a sharpened version of the Auslander-Reiten conjecture over $R\ltimes k$.    Also, we give a characterization of   the Betti numbers of an $R$-module over the idealization 
 ring $R\ltimes M$ and, as a biproduct,  we derive  that the Jorgensen-Leuschke conjecture holds true for $R\ltimes M$. Further, we show that the true of Buchsbaum-Eisenbud-Horrocks and  Total Rank conjectures over $R$ implies the true over $R\ltimes M$. This establishes particular answers for both conjectures for modules with infinite projective dimension, especially when $R$ is regular or a complete intersection ring.   As  applications of the idealization ring theory, we show that the  
 Zariski-Lipman conjecture holds  for any ring $R$ provided the Betti numbers  of the $R$-derivation module $\Der_k(R)$, seen as  $R\ltimes k$-module, satisfy the inequality
$\beta_{n}^{R\ltimes k}(\Der_k(R))\leq\beta_{n-1}^{R\ltimes k}(\Der_k(R))$  for some $n>0$. Some implications regarding the Herzog-Vasconcelos conjecture are also provided. 

\end{abstract}

\maketitle

\section{Introduction}

Let $R$ be a commutative Noetherian ring with identity and $M$ an $R$-module. The idealization ring $R\ltimes M$ of $M$ is the commutative ring $R\oplus M$ with identity $(1,0)$,  coordinate-wise addition and product defined by
$$(r_1,m_1)(r_2,m_2):=(r_1 r_2,r_1 m_2 +r_2 m_1).$$
This construction was introduced by Nagata (see \cite{Nagata}), and has been investigated by several authors  in different  subjects \cite{Winders, Gulliksen, Herzog, Nasseh-Yoshino, Reiten, Yamagishi}. Since we can put $M$ inside the commutative ring $R\ltimes M$ as an ideal, the theory of idealization of modules are useful for reducing results concerning submodules to the ideal case, generalizing facts from rings to modules and giving examples of commutative rings with zero divisors.

The main focus of this paper is to give positive answers for some  long-standing homological conjectures over the idealization ring $R\ltimes M$. In this direction, as the first open problem to be investigated in this work, we consider a general version of the celebrated Auslander-Reiten conjecture  \cite[p. 70]{AR}, posed in  \cite{TVRS}:

\begin{conjecture}[General Auslander-Reiten]\label{question1} For $M$ and $N$ finitely generated modules over a local ring $R$,
can one find integers $s$ and $t$, with $1\le s \le t$, such that for all $i$,
 $\Ext^i_R(M,N\oplus R)=0$   with $ s\leq i \leq t$, gives that 
  $M$ or $N$ has finite projective dimension?
\end{conjecture}

Moreover, this paper is particularly focused on investigating the ordinary version of the Auslander-Reiten conjecture:

\begin{conjecture}[Auslander-Reiten]\label{ARConjecture} Let $R$ be a commutative Noetherian local ring. If $M$ is a finitely generated $R$-module such that
$$\Ext_{R}^{i}(M,M\oplus R)=0$$
for all $i\geq 1$, then $M$ is free.
\end{conjecture}

Concerning the General Auslander-Reiten conjecture,  there are several situations in which the vanishing of
$\Ext^i_R(M,M\oplus R)$ for a specific finite set of values of $i$ is enough
 to deduce that $M$ is free or that $M$ 
has finite projective dimension (for instance \cite[Main Theorem]{Huneke-Leusch}, \cite[Corollary 10]{Araya} and \cite[Theorem 1.5]{Goto-Tak}). In the case that $M\neq N$ some answers for Conjecture \ref{question1} are given in \cite{TVRS}. For the ordinary case, it is well known in the literature that the Auslander-Reiten conjecture have several answers (see for instance   \cite{AINS,Celikbas-Tak,Christensen,K-O-T,NS}). In the idealization ring theory, Nasseh and Yoshino  \cite[Corollary 3.6]{Nasseh-Yoshino}  showed that the Auslander-Reiten conjecture is true over  $R\ltimes k$   provided $R$ is an Artinian local ring. Our  first main result, Theorem \ref{freenessmod} gives a positive answer for  Question \ref{question1} assuming an
additional Tor-rigid hypothesis regarding one of the modules. In addition, in Theorem \ref{ARC} we improve the result of Nasseh and Yoshino \cite[Corollary 3.6]{Nasseh-Yoshino} using a different approach,  removing the assumption on the ring $R$, assuming the vanishing of $\Ext^i_{R\ltimes k}(M,M\oplus (R\ltimes k))$ for some $i\geq 3$ and thus showing the true of the Auslander-Reiten Conjecture  for the idealization ring $R\ltimes k$.

In order to investigate  other   homological conjectures, first we give a characterization of the Betti numbers of any $R$-module over the idealization ring $R\ltimes M$ (Proposition \ref{bettinumbers=}). As a consequence, a refinement on the structure of  $R\ltimes M$  is furnished (see Theorem \ref{theostructure}, Proposition \ref{CIformula} and Corollary \ref{CIidealization}). As the main result, we cite that the idealization ring $R\ltimes M$ is always singular (Theorem \ref{theostructure} (i)). Also, we investigate the following  open problem posed by Jorgensen and Leuschke \cite[Question 2.6]{Jor-Leusch}.

\begin{conjecture}[Jorgensen-Leuschke]\label{JLquestion} Let $R$ be a Cohen-Macaulay local ring with canonical module $\omega$. Does $\beta_{1}^{R}(\omega)\leq\beta_{0}^{R}(\omega)$ imply that $R$ is Gorenstein?
\end{conjecture}

Some positive answers for this interesting conjecture  are given in  \cite{Cooper-Sather, paperbetti, Jor-Leusch, Lyle-Montano}. In this work we contribute to this investigation, showing that the Jorgensen-Leuschke Conjecture holds for any idealization ring $R\ltimes M$ (Theorem \ref{Jorleusch}). Regarding the behavior of the Betti numbers, other two famous conjectures can be considered:

\begin{conjecture}[Buchsbaum-Eisenbud-Horrocks] Let $R$ be a $d$-dimensional Noetherian local ring, and let $M$ be a finitely generated nonzero $R$-module. If $M$ has finite 
length and finite projective dimension, then for each $i \geq 0$ the $i$-th Betti number of $M$ over
$R$ satisfies the inequality
$$\beta_{i}^{R}(M)\geq\binom{d}{i}.$$
\end{conjecture}

 This famous open problem has few positive answers in the literature (\cite{Avramov2, Chang, Chara, Evans, paperbetti, Santoni}). In this sense, Avramov and Buchweitz \cite{Avramov2} introduced a weaker version of the previous conjecture, called the Total Rank Conjecture  as follows:
\begin{conjecture}[Total Rank] Let $R$ be a $d$-dimensional Noetherian local ring, and let $M$ be a finitely generated nonzero R-module. If $M$ has finite
length and finite projective dimension, then
$$\sum_{i\geq 0}\beta_{i}^{R}(M)\geq2^{d}.$$
\end{conjecture}

In \cite{Avramov2}, Avramov and Buchweitz showed that this conjecture  holds  for local rings of dimension 5 that contain their residue field. Also,  Walker  \cite{Walker} recently  has shown the true of Total rank Conjecture for  complete intersection rings $R$ such that  $2$ is invertible in $R$, or $R$ contains $\mathbb{Z}/p\mathbb{Z}$ as a subring for an odd prime $p$.  It is important to realize that these both conjectures  can be considered when the module does not have finite projective dimension. The known results in this topic were given by Tate \cite{tate} for hypersurfaces, and by Freitas-Jorge P\'erez for fiber product rings \cite{paperbetti}. In our context, we derive  that the true of the Buchsbaum-Eisenbud-Horrocks and  Total Rank conjectures over $R$ implies the true over $R\ltimes M$ (Theorem  \ref{BEH}). This provides particular answers for both conjectures for modules with infinite projective dimension, especially when $R$ is regular or a complete intersection ring.

In the last part  of this paper, we restrict our investigation to the derivation modules $\Der^p_S(R)$, where  $p\geq 1$  is an integer and $S$ is a commutative ring. The study  of derivation modules and their connections with the regularity of a variety have been motived by the famous Zariski-Lipman conjecture \cite{Lipman}, which states the following: Let \( X \) be a complex variety such that the tangent sheaf \( \mathcal{T}_X :=\mathcal{H}om_{\mathcal{O}_X}\left(\Omega_{X}^{(1)},\mathcal{O}_X\right) \) is locally free. Then \( X \) is smooth or non-singular. This conjecture remains widely open, but has been resolved in many special cases by various authors. References to the Zariski-Lipman Conjecture can be found in, e.g., Miller-Vassiliadou \cite[Section 4]{Milher}, as well as in the works of K\"allstr\"om \cite{Kallstrom} and Biswas, Gurjar, and Kolte \cite{Biswas}. Furthermore, this conjecture was stated in a more general form for $k$-algebras as follows: \begin{conjecture}[Zariski-Lipman] \label{conjecture2} Let $R$ be a Noetherian local ring. Then $R$ is a regular local ring if and only if $\Der_k(R)$ is a free $R$-module.
\end{conjecture}

For higher-order derivation modules,  Ludington \cite{Lu} posed the natural question:

\begin{question} \label{conjecture32} (Generalized Zariski-Lipman Conjecture) Assume that for some integer \( p\geq 1 \), \( \operatorname{pd}_{R} \operatorname{Der}_S^p(R) < \infty \). Under what assumptions on \( R \) and $S$, and for which values of \( p \) does this imply that \( R \) is regular? 
\end{question}

Our main result (Theorem \ref{mainteo}) gives a numerical criterion for the true of Question \ref{conjecture32} for any ring $R$.  In particular, we derive that the Zariski-Lipman conjecture holds true for any ring $R$ provided 
$\beta_{n}^{R\ltimes k}(\Der_k(R))\leq\beta_{n-1}^{R\ltimes k}(\Der_k(R))$ for some $n>0$ (Corollary \ref{ZLCclassica}). This result  strongly shows how the idealization ring theory can be a useful tool in obtaining answers to homological conjectures in general.  Following the same direction as the Zariski-Lipman Conjecture, there is also a homological conjecture proposed independently by Herzog and Vasconcelos: If ${\rm pd}_R{\rm Der}_k(R) \, < \, \infty$ then ${\rm Der}_k(R)$ is free. Differently to the Zariski-Lipman conjecture, the Herzog-Vasconcelos conjecture seems to be widely open, with some exceptions in specific cases; see, for example, \cite[Section 4]{He-M} and \cite{HV,CP, CT}. Our obtained results give some implications regarding this conjecture (see Corollaries \ref{HVconj} and \ref{HVconj2} for details).

The paper is organized as follows:  Section \ref{section2} establishes  our terminology and  some known results that will be utilized throughout the paper. Section \ref{section3} gives a criterion for the freeness of modules via vanishing of $\Ext$, in particular, answers for Questions \ref{question1} and \ref{ARConjecture}.  Section \ref{section4} deals with the Betti numbers of any $R$-module over the idealization ring $R\ltimes M$, and  in addition, we obtain conditions of when this rings are regular, hypersurfaces and complete intersection. Section \ref{section5} provides answers for the Jorgensen-Leuschke, Buchsbaum-Eisenbud-Horrocks and Total rank conjectures. Section \ref{section6} investigates the Generalized Zariski-Lipman conjecture over any ring $R$, using the idealization ring theory. Some consequences concerning the Herzog-Vasconcelos conjecture are also given.

\section{Preliminaries}\label{section2}

In this section, we review essential definitions and  results for the development of the paper. First, we list some key properties and facts of the theory of idealization rings that can be found in \cite{Winders}.
\begin{remark}\label{setupremark}
\begin{itemize}
    \item[(i)] If $(R,\m,k)$ is a local ring and $M$ is a finitely generated $R$-module, the idealization $R\ltimes M$ is a local ring with maximal ideal $\m\ltimes M$ and residue field $(R\ltimes M)/(\m\ltimes M)\cong R/\m=k$. 
    \item[(ii)] The Krull dimension of the idealization ringis equal the  Krull dimension of $R$, that is, $\dim(R\ltimes M)=\dim(R)$.
    \item[(iii)] If $R$ is Noetherian and $M$ is finitely generated, then $R\ltimes M$ is Noetherian.
    \item[(iv)] With the last conditions over $R$ and $M$, one has
    $$\depth_{R\ltimes M}(R\ltimes M)=\min\{\depth(R),\depth_{R}(M)\}$$
    and consequently, $R\ltimes M$ is Cohen-Macaulay if and only if $R$ is Cohen-Macaulay and $\depth_{R}(M)\geq\depth_{R}(R)$.
    \item[(v)] The $\m\ltimes M$-adic completion $\widehat{R\ltimes M}$ of $R\ltimes M$ is naturally isomorphic to $\widehat{R}\ltimes\widehat{M}$, where $\widehat{R}$ and $\widehat{M}$ are the $\m$-adic completion of $R$ and $\m$-adic completion of $M$, respectively.
    \item[(vi)] If $M=0$, $R\ltimes M$ is naturally isomorphic to $R$ and if $M=R^{n}$ is a free $R$-module, then $R\ltimes M$ is  isomorphic to $R[X_1,\dots,X_n]/(X_1,\dots,X_n)^{2}$.

    \item[(vii)]  Considering the natural inclusion and projection maps, there is a exact sequence of $R\ltimes M$-modules given by
 $0\xrightarrow[]{}M\xrightarrow[]{}R\ltimes M\xrightarrow[]{}R\xrightarrow[]{}0.
$
\end{itemize}
\end{remark}

\begin{remark}\label{ringhomomorphism} We can consider an $R\ltimes M$-module as an $R$-module via the homomorphism of rings $R\xrightarrow[]{}R\ltimes M$ given by $r\to (r,0)$. Note that the composition of this homomorphism with the projection $R\ltimes M\xrightarrow[]{}R$ is the identity mapping on $R$. We utilize these local homomorphisms to establish that if $L$ is an $R$-module (resp. an $R \ltimes M$-module), then $L$ is also an $R \ltimes M$-module (resp. an $R$-module), respectively.
\end{remark}

\begin{notation}{\rm 

For the rest of the paper, the ring $R$ is a commutative  Noetherian local ring with maximal ideal $\m$, residue field $k$ and the $R$-modules are nonzero and finitely generated. Also, all $R\ltimes M$-modules are nonzero and finitely generated. Let $\mu_{R}(M)$ denote the minimal number of generators of an $R$-module $M$. The embedding dimension of  $R$, denoted by $\operatorname{edim}(R)$, is the minimal number
of generators of its maximal ideal $\m$, i.e., $\operatorname{edim}(R):=\dim_{k}(\m/\m^{2})$. In general, $\operatorname{edim}(R)\geq\dim(R)$ and the equality happens provided $R$ is a regular ring.
}
\end{notation}

\begin{definition} Let $R$ be a ring and let $M$ be an $R$-module. The formal power series
$$P_{M}^{R}(t):=\sum_{i\geq 0}\beta_{i}^{R}(M)t^{i}$$
is called \textit{Poincar\'{e} series} of $M$, where $\beta_{i}^{R}(M):=\dim_{k}\Tor_{i}^{R}(M,k)$ is the $i$-th Betti number of $M$, that is, the dimension of the $k$-vector space $\Tor_{i}^{R}(M,k)$. 

\end{definition}



\section{Vanishing of  Ext over idealization rings }\label{section3}

In this section, the main focus is the investigation of the general version of the Auslander-Reiten conjecture described in the introduction (Question \ref{question1}). For this purpose,  consider the right projective approximation $P_{0}\xrightarrow[]{\partial_{0}}M$. We define the first \textit{syzygy} of $M$ as $\Omega_{R}^{1}M:=\ker\partial_{0}$, which is unique up to projective equivalence. Additionally, we can recursively define the $n$-th syzygy module of $M$ as $\Omega_{R}^{n}M:=\Omega_{R}^{1}(\Omega_{R}^{n-1}M)$ for all $n\geq 1$, with $\Omega_R^0M:=M$ by convention.

Let $M^{\ast}:=\Hom_{R}(M,R)$ be the dual of $M$. The {\em Auslander transpose}, written as $D M$ or $D_1M$,   is defined as the cokernel of the map $F_0^*\to F_1^*$, where 
$$
F_1\to F_0\to M\to 0
$$ is a minimal resolution of $M$, with  $F_i$ free $R$-modules.  Thus one has the exact sequence
\begin{equation}\label{Austranspose}
0 \to M^* \to F_0^*\to F_1^* \to D M \to 0\,.
\end{equation}

More generally, for  $n\geq 1$ and  a minimal free resolution $F$ of $M$  over~$R$,
$$
F: \cdots \to F_n\to\cdots  \to F_2\to F_1\to F_0\to M\to 0,
$$
define the $n^{\text {th}}$  Auslander transpose  $D_nM$ by  $D_nM:=\Coker\left(F_{n-1}^*\to F_n^*\right)$.


An important result related to these definitions is the existence of the following:

\begin{remark}[Auslander sequence](\cite[Theorem 2.8 (b)]{Auslander-Bridger} or \cite[(1.1.1)]{D-J}) For any finitely generated $R$-modules $M$ and $N$ and for all $i\geq 0$, there exists an exact sequence
\begin{align}\label{ausseq}\operatorname{Tor}_{2}^{R}(D\Omega_{R}^{i}M,N)&\xrightarrow[]{}\operatorname{Ext}_{R}^{i}(M,R)\otimes_{R}N &\\
    &\xrightarrow[]{}\operatorname{Ext}_{R}^{i}(M,N)\xrightarrow[]{}\operatorname{Tor}_{1}^{R}(D\Omega_{R}^{i}M,N)\xrightarrow[]{}0.\nonumber &
    \end{align}
\end{remark}

 By the canonical map $h:M\xrightarrow[]{}M^{\ast\ast}$, \cite[Exercise 1.4.21]{Bruns-Herzog} provides the exact sequence
\begin{equation}\label{reflexseq} 0\xrightarrow[]{}\Ext_{R}^{1}(DM,R)\xrightarrow[]{}M\xrightarrow[]{h}M^{\ast\ast}\xrightarrow[]{}\Ext_{R}^{2}(DM,R)\xrightarrow[]{}0.
\end{equation}
Note that $M$ is reflexive if and only if $\Ext_{R}^{1}(DM,R)=\Ext_{R}^{2}(DM,R)=0$.


Now, we start with the following key result given by Nasseh and Yoshino \cite{Nasseh-Yoshino}.

\begin{lemma}\cite[Theorem 3.1]{Nasseh-Yoshino}\label{N-YTor} Let $(R,\m,k)$ be a local ring and $M$ and $N$ be nonzero finitely generated $R\ltimes k$-modules such that $\Tor_{n}^{R\ltimes k}(M,N)=0$ for some $n\geq 3$. Then $M$ is $R\ltimes k$-free or $N$ is $R\ltimes k$-free.
\end{lemma}

\begin{definition}\cite{Auslander1} An $R$-module $M$ is said to be Tor-rigid provided  the following holds for every $R$-module $N$ and every $n\geq 1$
$$\Tor_{i}^{R}(M,N)=0\Longrightarrow\Tor_{i+1}^{R}(M,N)=0.$$
\end{definition}



We are able to show one of the main result of this section.

\begin{theorem}\label{freenessmod} Let $M$ and $N$ be  $R\ltimes k$-modules. If $N$ is Tor-rigid and $\Ext_{R\ltimes k}^{1}(M,N)=0$, then $M$ is free or $N$ is free.
\end{theorem}

\begin{proof} Suppose that $\Ext_{R\ltimes k}^{1}(M,N)=0$. This combined with the Auslander sequence \ref{ausseq} gives
$\Tor_{1}^{R\ltimes k}(D\Omega_{R\ltimes k}^{1}M,N)=0.$
Since $N$ is Tor-rigid one has
\begin{equation}\label{Toreq} \Tor_{i}^{R\ltimes k}(D\Omega_{R\ltimes k}^{1}M,N)=0 \text{ for all }i>0.    
\end{equation}
Then, by Lemma \ref{N-YTor} we have that $N$ is $R\ltimes k$-free or $D\Omega_{R\ltimes k}^{1}M$ is $R\ltimes k$-free. In the latter condition, $D\Omega_{R\ltimes k}^{1}M$ is $R\ltimes k$-free if and only if $\Omega_{R\ltimes k}^{1}M$ is $R\ltimes k$-free, and then $\pd_{R\ltimes k}(M)<\infty$. Therefore $M$ is $R\ltimes k$-free by the Auslander-Buchsbaum formula. 
\end{proof}

As a consequence, we derive a positive answer for  the Question \ref{question1} assuming an additional Tor-rigid hypothesis regarding one of the modules.

\begin{corollary} Let $M$ and $N$ be  $R\ltimes k$-modules. If $N$ is Tor-rigid and $\Ext_{R\ltimes k}^{i}(M,N)=0$ for some $i>0$, then $M$ is free or $N$ is free.
\end{corollary}

\begin{proof} Since $\Ext_{R\ltimes k}^{1}\left(\Omega_{R\ltimes k}^{i}M,N\right)=\Ext_{R\ltimes k}^{i}(M,N)=0$, Theorem \ref{freenessmod} gives that $\Omega_{R\ltimes k}^{i}M$ is free or $N$ is free. As the depth of $R\ltimes k$ is zero (Remark \ref{setupremark} (iv)), the freeness of $M$  or $N$  follows by the Auslander-Buchsbaum formula \cite[A.5. Theorem, p. 310]{LW}.
\end{proof}

The next result improves \cite[Corollary 3.5]{Nasseh-Yoshino}.

\begin{proposition}\label{NYExt} Let $M$ and $N$ be $R\ltimes k$-modules. Suppose that $\Ext_{R\ltimes k}^{i}(M,N)=0$ for some $i\geq\max\{3,\depth_{R\ltimes k}(N)+3\}$. Then $M$ is $R\ltimes k$-free or $N$ is $R\ltimes k$-injective.
\end{proposition}

\begin{proof} Suppose that $\Ext_{R\ltimes k}^{i}(M,N)=0$ for some $i\geq\max\{3,\depth_{R\ltimes k}(N)\}$ and assume that $M$ is not free. Since $\depth(R\ltimes k)=0$, then $M$ has infinite projective dimension. Note that $(0,1)\in (0:\m)$ but $(0,1)\notin\m^{2}$. By \cite[Proposition 3.3(2)]{AINS} we obtain that $k$ is a direct summand of $\Omega_{R\ltimes k}^{2}(M)$. Hence we obtain $\Ext_{R\ltimes k}^{i-2}(k,N)=0$. It follows by \cite[Theorem 2]{Roberts} that $\id_{R\ltimes k}(N)<\infty$. Now,  the Bass formula gives that $N$ is injective.
\end{proof}

In \cite[Corollary 3.6]{Nasseh-Yoshino}, Nasseh and Yoshino showed that the Auslander-Reiten conjecture holds true for the idealization ring $R\ltimes k$, provided $R$ is an Artinian local ring.  As the second  main result of this section, we  improve the result of Nasseh and Yoshino, and thus providing a positive answer for Question \ref{ARConjecture}.

\begin{theorem}\label{ARC} Let $M$ be a  $R\ltimes k$-module. If $\Ext_{R\ltimes k}^{i}(M,M\oplus (R\ltimes k))=0$ for some $i\geq 3$, then $M$ is free.
\end{theorem}

\begin{proof} Suppose that $\Ext_{R\ltimes k}^{i}(M, M \oplus (R\ltimes k))=0$ for some $i\geq 3$ and assume that $M$ is not free. With the same steps in the proof of Proposition \ref{NYExt}, we obtain that $k$ is a direct summand of $\Omega_{R\ltimes k}^{i-1}(M)$. So $\Ext_{R\ltimes k}^{1}(k,R\ltimes k)=0=\Ext_{R\ltimes k}^{1}(k,M)$ and $R\ltimes k$ is Gorenstein by \cite[Theorem 2]{Roberts}. Therefore, it follows by Proposition \ref{NYExt} that $M$ is $R\ltimes k$-injective. Therefore $M$ has finite projective dimension, a contradiction.
\end{proof}

\begin{remark}\label{fiberproduct} It is important to realize that the range of $i$-values in Lemma \ref{N-YTor} is optimal from \cite[Remark 3.4]{Nasseh-Yoshino}. In addition, since $R\ltimes k\cong R\times_{k}(k[x]/(x^{2}))$ is a fiber product ring of depth 0, the Auslander-Reiten conjecture holds true provided   $i\geq 0$ (\cite[Theorem 1.1]{NS}). Recently, in \cite{TVRS}, the authors generalize \cite[Theorem 1.1]{NS} and  showed the true of Auslander-Reiten conjecture for fiber product rings for $1\leq i \leq 6$.   

\end{remark}

The last part of this section is devoted to obtain a relation between reflexitivity and freeness of modules using the vanishing of certain Ext modules.

\begin{lemma}\label{reflexive} Let  $M$ be an $R$-module. If $M$ is reflexive and $\Ext_{R}^{1}(M,R)=0$, then $\Omega_{R}^{1}M$ is reflexive. In particular, if $M$ is reflexive and $\Ext_{R}^{i}(M,R)=0$ for  $1\leq i\leq n$, then $\Omega_{R}^{i}M$ is reflexive for  $1\leq i\leq n$.
\end{lemma}

\begin{proof} Suppose that $M$ is reflexive. From \ref{reflexseq} we obtain $\Ext_{R}^{2}(DM,R)=0$. Applying this in the exact sequence
$$0\xrightarrow[]{}(\Omega_{R}^{1}M)^{\ast}\xrightarrow[]{}F_{1}\xrightarrow[]{}DM\xrightarrow[]{}0,$$
we derive
\begin{equation}\label{syzext} \Ext_{R}^{1}((\Omega_{R}^{1}M)^{\ast},R)=0.
\end{equation}
Since $\Ext_{R}^{1}(M,R)=0$, one has another exact sequence
$$0\xrightarrow[]{}M^{\ast}\xrightarrow[]{}F_{0}^{\ast}\xrightarrow[]{}(\Omega^{1}M)^{\ast}\xrightarrow[]{}0,$$
and combining this with \eqref{syzext} provides
$$0\xrightarrow[]{}(\Omega_{R}^{1}M)^{\ast\ast}\xrightarrow[]{}F_{1}^{\ast\ast}\xrightarrow[]{}M^{\ast\ast}\xrightarrow[]{}0.$$
Therefore, the five's lemma applied on the diagram
\begin{center}\begin{tikzcd} 0\arrow{r} & \Omega_{R}^{1}M\arrow{r}\arrow{d} & F_{1}\arrow{r}\arrow{d} & M\arrow{r}\arrow{d} & 0 &\\
0\arrow{r} & (\Omega_{R}^{1}M)^{\ast\ast}\arrow{r} & F_{1}^{\ast\ast}\arrow{r} & M^{\ast\ast}\arrow{r} & 0
\end{tikzcd}
\end{center}
implies that $\Omega_{R}^{1}M$ is reflexive.
\end{proof}


\begin{proposition} Let $M$ be an $R\ltimes k$-module. Then, $M$ is reflexive and $\Ext_{R\ltimes k}^{1}(M,R\ltimes k)=0$ if and only if $M$ is free.
\end{proposition}

\begin{proof} Suppose that $M$ is reflexive. From \ref{reflexseq} one has  $\Ext_{R\ltimes k}^{2}(DM,R\ltimes k)=0$. By \eqref{Austranspose} $\Omega_{R\ltimes k}^{2}D\Omega_{R\ltimes k}^{1}M=(\Omega_{R\ltimes k}^{1}M)^{\ast}$ and by construction and the fact that $\Ext_{R\ltimes k}^{1}(M,R\ltimes k)=0$, we have $\Omega_{R\ltimes k}^{1}DM=(\Omega_{R\ltimes k}^{1}M)^{\ast}$. Hence
\begin{align*} 0 &= \Ext_{R\ltimes k}^{2}(DM,R\ltimes k)=\Ext_{R\ltimes k}^{1}(\Omega_{R\ltimes k}^{1}DM,R\ltimes k) &\\
&= \Ext_{R\ltimes k}^{1}(\Omega_{R\ltimes k}^{2}D\Omega_{R\ltimes k}^{1}M,R\ltimes k)=\Ext_{R\ltimes k}^{3}(D\Omega_{R\ltimes k}^{1}M,R\ltimes k).
\end{align*}
Consider $N:=D\Omega_{R\ltimes k}^{1}M$ and suppose that $N$ is not free. As $\depth(R\ltimes k)=0$, we have that $\pd_{R\ltimes k}(N)=\infty$. It follows by \cite[Proposition 3.3(2)]{AINS} that $k$ is a direct summand of $\Omega_{R\ltimes k}^{2}N$. Then
$$\Ext_{R\ltimes k}^{1}(\Omega_{R\ltimes k}^{2}N,R\ltimes k)=\Ext_{R\ltimes k}^{3}(N,R\ltimes k)=0$$
and this implies that $\Ext_{R\ltimes k}^{1}(k,R\ltimes k)=0$. Now, \cite[Theorem 2]{Roberts} gives  that $\id_{R\ltimes k}(R\ltimes k)=0$, forcing $R\ltimes k$ to be a Gorenstein ring. However, by Remark \ref{fiberproduct} and \cite[Corollary 2.7]{NSTV} $R\ltimes k$ can not be a Gorenstein ring. Then $D\Omega_{R\ltimes k}^{1}M$ is free, implying that $(\Omega_{R\ltimes k}^{1}M)^{\ast}$ is free. Therefore, Lemma \ref{reflexive} provides that  $\Omega_{R\ltimes k}^{1}M$ is reflexive and, consequently, $\Omega_{R\ltimes k}^{1}M=(\Omega_{R\ltimes k}^{1}M)^{\ast\ast}$ is free. The Auslander-Buchsbaum formula \cite[A.5. Theorem, p. 310]{LW} gives that  $M$ is free, because the depth of $R\ltimes k$ is zero.
\end{proof}


\section{Betti Numbers and Structural results of $R\ltimes M$ }\label{section4}


An important notion for the rest of this work is the characterization of the Poincar\'e series of any $R$-module over the idealization ring $R\ltimes M$ given by Herzog \cite{Herzog}.




\begin{lemma}\cite[Corollary 1]{Herzog}\label{poincareidealization} Let $R$ be a local ring and let $M$ and $N$ be $R$-modules. Then
$$P_{N}^{R\ltimes M}(t)=\frac{P_{N}^{R}(t)}{1-tP_{M}^{R}(t)}.$$
\end{lemma}

The next result provides  the  shape of the Betti numbers of an $R$-module $N$ over $R\ltimes M$. This characterization is a key ingredient for the investigation of the structural results and homological conjectures given below.
\begin{proposition}\label{bettinumbers=} Let  $M$ and $N$ be $R$-modules. Then
\begin{itemize}
    \item[(i)] $\beta_{0}^{R\ltimes M}(N)= \beta_{0}^{R}(N)$.
    \item[(ii)] $\beta_{1}^{R\ltimes M}(N)= \beta_{0}^{R}(N)\beta_{0}^{R}(M)+\beta_{1}^{R}(N)$. In particular, the following equality happens  $\operatorname{edim}(R\ltimes M)= \beta_{0}^{R}(M)+\operatorname{edim}(R)$.
    \item[(iii)] $\beta_{2}^{R\ltimes M}(N)= \beta_{0}^{R}(N)\beta_{0}^{R}(M)^{2}+\beta_{0}^{R}(N)\beta_{1}^{R}(M)+\beta_{0}^{R}(M)\beta_{1}^{R}(N)+\beta_{2}^{R}(N)$.
    \item[(v)] $\beta_{i}^{R\ltimes M}(N)\geq \beta_{0}^{R}(M)\beta_{i-1}^{R}(N)+\beta_{i}^{R}(N)$ for all $i\geq 0$.
   
\end{itemize}
\end{proposition}

\begin{proof} By Theorem \ref{poincareidealization}, the Poincar\'{e} series $P_{N}^{R\ltimes M}(t)$ is given by
$$P_{N}^{R\ltimes M}(t)=\frac{\sum_{i\geq 0}\beta_{i}^{R}(N)t^{i}}{1-\sum_{i\geq 1}\beta_{i-1}^{R}(M)t^{i}}.$$
Let
$$\frac{1}{1-\sum_{i\geq 1}\beta_{i-1}^{R}(M)t^{i}}:=\sum_{i\geq 0}B_{i}t^{i}.$$
Let $b_{i}$ be the coeficients of $1-\sum_{i\geq 1}\beta_{i-1}^{R}(M)t^{i}$. Note that $b_{0}=1$ and $b_{i}=-\beta_{i-1}^{R}(M)$ for $i\geq 1$. Since $b_{0}\neq 0$, we have $B_{0}=1$ and for all $i\geq 1$, $B_{i}$ is given by the determinant
$$B_{i}=\frac{1}{i!}\begin{vmatrix}
0 & ib_{1} & ib_{2} & ib_{3} & \cdots & ib_{i} \\
0 & (i-1) & (i-1)b_{1} & (i-1)b_{2} & \cdots & (i-1)b_{i-1} \\
0 & 0 & (i-2) & (i-2)b_{1} & \cdots & (i-2)b_{i-2} \\
0 & 0 & 0 & (i-3) & \cdots & (i-3)b_{i-3} \\
\vdots & \vdots & \vdots & \vdots & \ddots & \vdots \\
1 & 0 & 0 & 0 & \cdots & 1 \\
\end{vmatrix}.$$
Moreover, $B_{i}=\sum_{j=1}^{i}\left| b_{j} \right|B_{i-j}$. Since $b_{i}\leq 0$ for $i\geq 1$, then $B_{i}\geq 0$ for all $i\geq 0$. Hence
$$P_{N}^{R\ltimes M}(t)=\left ( \sum_{i\geq 0}\beta_{i}^{R}(N) \right )\left ( \sum_{i\geq 0}B_{i}t^{i} \right )=\sum_{n\geq 0}\left ( \sum_{i=0}^{n}\beta_{i}^{R}(N)B_{n-i} \right )t^{n},$$
which allows us to conclude that
\begin{equation}\label{betti} \beta_{n}^{R\ltimes M}(N)=\sum_{i=0}^{n}\beta_{i}^{R}(N)B_{n-i}.    
\end{equation}
The statements $(i)-(iv)$ follows due to equation \ref{betti}.
\end{proof}

\begin{remark}\label{pdinifinitaideal}
Note that the previous result furnishes that even if $\pd_{R}(N)\leq\infty$ is not true that $\pd_{R\ltimes M}(N)<\infty$. Considering the case where $N$ is a free $R$-module, since $\beta_{i}^{R\ltimes M}(N)\geq \beta_{0}^{R}(M)\geq 1$ for $i\geq 1$, we see that $\pd_{R\ltimes M}(N)=\infty$.  
\end{remark}

Now, we are able to  study the structure of the idealization ring $R\ltimes M$. 
The Cohen-Macaulayness, and the Gorensteiness  of  $R\ltimes M$ were given in \cite{Winders} and \cite{Reiten}. As a refinement of this invesgation, we give the following results.

\begin{theorem}\label{theostructure} Let $M$ be an $R$-module. Then the following happens.
\begin{itemize}
    \item[(i)]  $R\ltimes M$ is not regular.
    
    \item[(ii)] If $R\ltimes M$ is a hypersurface, then $M$ is a cyclic $R$-module and $R$ is regular. The converse holds provided $R\ltimes M$ is Cohen-Macaulay.

\end{itemize}
    
\end{theorem}

\begin{proof} First, note that $\dim(R\ltimes M)=\dim(R)$ (Remark (ii) \ref{setupremark}).

(i) Assume that $R\ltimes M$ is regular. By Proposition \ref{bettinumbers=} (ii)
$$\dim(R)=\dim(R\ltimes M)=\operatorname{edim}(R\ltimes M)=\beta_{0}^{R}(M)+\operatorname{edim}(R).$$
Since $\operatorname{edim}(R)\geq \dim(R)$, one has that $\beta_{0}^{R}(M)=0$, and so $M=0$, a contradiction. Therefore $R\ltimes M$ is not regular.

(ii) Suppose that $R\ltimes M$ is a hypersurface. Hence $\operatorname{edim}(R\ltimes M)-\depth(R\ltimes M)\leq 1$ and since $R\ltimes M$ is Cohen-Macaulay one obtains $\depth(R\ltimes M)=\dim(R\ltimes M)=\dim(R)$. By Proposition \ref{bettinumbers=} (ii)
$$\beta_{0}^{R}(M)+\operatorname{edim}(R)-\dim(R)=\operatorname{edim}(R\ltimes M)-\depth(R\ltimes M)\leq 1.$$
Then $\beta_{0}^{R}(M)=1$ and $\operatorname{edim}(R)=\dim(R)$ and consequently $R$ is regular.

On the other hand, since $R\ltimes M$ is Cohen-Macaulay, Proposition \ref{bettinumbers=} (ii) together with the hypothesis imply
$$1=\beta_{0}^{R}(M)+\operatorname{edim}(R)-\dim(R)=\operatorname{edim}(R\ltimes M)-\depth(R\ltimes M).$$
Therefore, $R\ltimes M$ is a hypersurface.
\end{proof}

\begin{proposition}\label{CIformula} Let  $M$ be an $R$-module. Then $R\ltimes M$ is a complete intersection provided
$$\frac{\beta_{0}^{R}(M)(1-\beta_{0}^{R}(M))+\beta_{1}^{R}(M)(1+\beta_{1}^{R}(M))}{\beta_{1}^{R}(M)+\beta_{2}^{R}(k)+\dim(R)}=2.$$
\end{proposition}

\begin{proof} First, note that $R$ and $R\ltimes M$ have the same residue field $k$ (Remark \ref{setupremark} (i)). By \cite[Theorem 7.3.3]{Avramov} (or \cite[Proposition 2.8.4 (3)]{Hashimoto}, \cite[Theorem 2.3.3]{Bruns-Herzog}), $R\ltimes M$ is a complete intersection if and only if
\begin{equation}\label{eq1} \beta_{2}^{R\ltimes M}(k)=\binom{\beta_{1}^{R\ltimes M}(k)}{2}+\beta_{1}^{R\ltimes M}(k)-\dim(R\ltimes M).   
\end{equation}
Now, Proposition \ref{bettinumbers=} (iii) gives that
\begin{equation}\label{eq2} \beta_{2}^{R\ltimes M}(k)=\beta_{0}^{R}(M)^{2}+\beta_{1}^{R}(M)+\beta_{0}^{R}(M)\beta_{1}^{R}(k)+\beta_{2}^{R}(k).
\end{equation}
Comparing \eqref{eq1} and \eqref{eq2}, the desired conclusion follows.
\end{proof}

The previous result shows that is not easy to obtain idealization rings $R\ltimes M$ that are complete intersection. However in the case that $R$ is complete intersection we derive the following consequence. 

\begin{corollary}\label{CIidealization} Let $M$ be a module over a complete intersection $R$. If $R\ltimes M$ is a complete intersection, then $M=R$. The converse holds provided $R$ is regular.
\end{corollary}

\begin{proof} First, note that if $R$ is a regular ring and $M=R$, it follows by Remark \ref{setupremark} (vi) that $R\ltimes R\cong R[X]/(X^{2})$, which is a complete intersection ring. 

Now, suppose that $R\ltimes M$ is a complete intersection. Since $R$ is a complete intersection, then  $\beta_{2}^{R}(k)=\binom{\beta_{1}^{R}(k)}{2}+\beta_{1}^{R}(k)-\dim(R)$ by \cite[Theorem 7.3.3]{Avramov} (or \cite[Proposition 2.8.4 (3)]{Hashimoto}, \cite[Theorem 2.3.3]{Bruns-Herzog}). Then, substituing this in the formula of Proposition \ref{CIformula}, we obtain that $R\ltimes M$ is a complete intersection provided
$$\beta_{1}^{R}(M)=\frac{\beta_{0}^{R}(M)-\beta_{0}^{R}(M)^{2}}{2},$$
and this happens if and only if $\beta_{0}^{R}(M)=1$ and $\beta_{1}^{R}(M)=0$, that is, $M=R$.
\end{proof}


\section{Homological Conjectures over idealization rings}\label{section5}

The aims of this section is to provide answers to some homological conjectures described in the introduction. First, we investigate the following  conjecture posed by Jorgensen and Leuschke \cite[Question 2.6]{Jor-Leusch}.

\begin{conjecture}    
[Jorgensen-Leuschke] Let $R$ be a Cohen-Macaulay ring with canonical module $\omega$. Does $\beta_{1}^{R}(\omega)\leq\beta_{0}^{R}(\omega)$ imply that $R$ is Gorenstein?
\end{conjecture}

In order to study  the previous open problem in the context of idealizations rings, it is natural to ask when $R\ltimes M$ is Gorenstein.    In \cite[Theorem 7]{Reiten}, Reiten showed that $R\ltimes M$ is a Gorenstein ring if and only if $R$ is a Cohen-Macaulay ring and $M$ is the canonical module of $R$. 
This key characterization helps us to show in the next result that the Jorgensen-Leuschke Conjecture is true for the idealization ring  $R\ltimes M$.







\begin{theorem}\label{Jorleusch} Let $R\ltimes M$ be Cohen-Macaulay ring with canonical module $\omega$. 
Suppose that  $\beta_{1}^{R\ltimes M}(\omega)\leq\beta_{0}^{R\ltimes M}(\omega)$.  Then $R\ltimes M$ is Gorenstein.
\end{theorem}

\begin{proof}

First,  we may consider $\omega$ as an $R$-module by  Remark \ref{ringhomomorphism}. In addition, since $\beta_{1}^{R\ltimes M}(\omega)\leq\beta_{0}^{R\ltimes M}(\omega)$, Proposition \ref{bettinumbers=} gives
$$\beta_{0}^{R}(\omega)=\beta_{0}^{R\ltimes M}(\omega)\geq\beta_{1}^{R\ltimes M}(\omega)=\beta_{0}^{R}(\omega)\beta_{0}^{R}(M)+\beta_{1}^{R}(\omega).$$
This yields $\beta_{1}^{R}(\omega)=0$ and $\beta_{0}^{R}(M)=1$, and so  $\omega$ is a free $R$-module. Since the Betti numbers are invariants by completion, we may assume that $R$ and $R\ltimes M$ are complete. By \cite[Corollary 5.14]{Herzog-Kunz} and the isomorphism in the proof of \cite[Theorem 3.3]{Yamagishi}, we derive the $R$-isomorphism
    \begin{equation}\label{isocanmod} \omega\cong\Hom_{R}(R\ltimes M,\omega_{R})\cong\Hom_{R}(M,\omega_{R})\oplus\omega_{R},
    \end{equation}
    where $\omega_{R}$ is the canonical module of $R$. Hence $\omega_{R}$ is a free $R$-module and therefore $R$ is a Gorenstein ring. In particular, we obtain that $\omega_{R}=R$ and the isomorphism \ref{isocanmod} furnishes that $M^{\ast}= \Hom_{R}(M,R)= \Hom_{R}(M,\omega_{R})$ is a free $R$-module.
    
Since $R\ltimes M$ is Cohen-Macaulay, \cite[Corollary 4.14]{Winders} provides that $R$ is Cohen-Macaulay and ${\rm depth}_R(M) \geq {\rm depth}(R)= \dim (R)$. Also $\dim (R)\geq {\rm depth}_R(M)$, and so $\dim (R)= {\rm depth}_R(M)$, i.e., $M$ is maximal Cohen-Macaulay $R$-module. Further,    \cite[Theorems 3.3.7 and 3.3.10(d)]{Bruns-Herzog} yields that $M$ is reflexive, because $R$ is Gorenstein.  The reflexivity of $M$ and the fact that $M^{\ast}$ is free  imply  the freeness of $M$.
Now, since $\beta_{0}^{R}(M)=1$, one obtains that  $M=R=\omega_{R}$. Therefore  $R\ltimes M$ is Gorenstein by \cite[Theorem 7]{Reiten}, as desired.  
 \end{proof}



Another on going topic of investigation in homological algebra are the lower bounds for the Betti numbers. In this sense, as described in the introduction, two long-standing open problems  are highlighted, the Buchsbaum-Eisenbud-Horrocks  and Total Rank  conjectures. From the characterization of the Betti numbers given in Proposition \ref{bettinumbers=}, the next result provides a positive answer for both conjectures over the idealization ring $R\ltimes M$ for modules that have infinite projective dimension (Remark \ref{pdinifinitaideal}).

\begin{theorem} \label{BEH} Let $R$ be a $d$-dimensional ring and let $M$, $N$ be $R$-modules.
\begin{itemize}
    \item[(i)] If $R$ satisfies the Buchsbaum-Eisenbud-Horrocks Conjecture and $N$ has finite lenght and finite projective dimension over $R$, then for all $1\leq i\leq \depth(R)$    $$\beta_{i}^{R\ltimes M}(N)\geq \binom{d}{i}.$$
    \item[(ii)] If $R$ satisfies the Total Rank Conjecture and $N$ has finite lenght and finite projective dimension over $R$, then for all $1\leq i\leq \depth(R)$    $$\sum_{i\geq 0}\beta_{i}^{R\ltimes M}(N)\geq 2^{d}.$$
\end{itemize}
\end{theorem}

\begin{proof} By Proposition \ref{bettinumbers=} (v) we derive that $\beta_{i}^{R\ltimes M}(N)\geq\beta_{i}^{R}(N)$. Therefore, since  $R$ satisfies the Buchsbaum-Eisenbud-Horrocks Conjecture, it follows that
$$\beta_{i}^{R\ltimes M}(N)\geq\beta_{i}^{R}(N)\geq\binom{d}{i}.$$
The proof of (ii) is similar to (i).
\end{proof}

It is well known in the literature that the Buchsbaum-Eisenbud-Horrocks  conjecture holds  for the residue field $k$ over a regular ring  $R$, and that any ring that satisfies the Buchsbaum-Eisenbud-Horrocks conjecture also satisfies the Total rank conjecture. In addition, Walker (\cite[Theorem 2 (1)]{Walker}) has shown the positivity of Total rank Conjecture  for  complete intersection rings where the characteristic  of the residue field $k$ is not two. With this facts, we provide the following consequences.
\begin{corollary}
      Let $R$ be a $d$-dimensional ring and  let $M$ be an $R$-module.
\begin{enumerate}

\item[(i)] If $R$ is regular, then for all $1\leq i \leq {\rm depth} (R)$
    $$\beta_i^{R\ltimes M}(k)\geq \binom{d}{i}.$$

    \item[(ii)] If $R$ is regular, then  $$\sum_{i= 0}^{ {\rm depth}(R)}\beta_i^{R\ltimes M}(k)\geq 2^{d}.$$

    \item[(iii)] Suppose that $R$ is the quotient of a regular local ring by a regular sequence of elements and the  characteristic of $k$ is different from 2. If $N$ is an $R$-module with finite projective dimension and finite length over $R$, then 
$$\sum_{i=0}^d\beta_i^{R\ltimes M}(N)\geq 2^d.$$
    
\end{enumerate}
\end{corollary}

\begin{remark}
Note that the previous result yields new examples of the positivity of the Buchsbaum-Eisenbud-Horrocks  and Total Rank conjectures, because  even though $R$ is regular or a complete intersection, the idealization ring $R\ltimes M$ is not regular (Theorem \ref{theostructure}) or eventually   a complete intersection ring. For instace, in the case that $M$ is not maximal Cohen-Macaulay, Remark \ref{setupremark} (iv) gives that $R\ltimes M$ is not Cohen-Macaulay, consequently, is not a complete intersection (see  Proposition \ref{CIformula}  for further details).       
\end{remark}

\section{ Zariski-Lipman and Herzog-Vasconcelos Conjectures}\label{section6} 
In this section, we give a brief addendum on the Zariski-Lipman and Herzog-Vasconcelos Conjecture. First, we fix our notation and recall some known results of  K\"ahler  differential modules and differential modules in the context of rings. 

\subsection*{Differential modules of order $p$.}  All the rings here are commmutative, Noetherian with identity. 
 Let $R$ be a $S$-algebra and $M$ be an $R$-module. Recall that a $S$-linear map $D : R \rightarrow M$ is said to be a derivation if for any two elements $x_0, x_1 \in R$, the following identity holds:
\[ D(x_0 x_1) = x_0 D(x_1) + x_1 D(x_0). \]
A derivation of order $n$ can be defined generalizing the previous identity as follows.

 A $S$-linear map $D : R \rightarrow M$ is said to be a $p$-{\it th order derivation} if for any $x_0, \ldots, x_n \in R$, the following identity holds:
\[ D(x_0 \cdots x_p) = \sum_{s=1}^{p} (-1)^{s-1} \sum_{i_1 < \cdots < i_s} x_{i_1} \cdots x_{i_s} D(x_0 \cdots \hat{x}_{i_1} \cdots \hat{x}_{i_s} \cdots x_p), \]
where $\hat{x}_{i_j}$ means that this element does not appear in the product.  The set of $n$-th order derivations of an $S$-algebra $R$ into an $R$-module $M$ over $S$ will be denoted by $\text{Der}^p_S(R, M)$. When $M = R$, we shall use the notation $\text{Der}^p_S(R)$ in place of $\text{Der}^p_S(R, R)$.

The module of derivations of order $p$ can be represented as follows. Let $I$ denote the kernel of the homomorphism $R \otimes_S R \rightarrow R$, $a \otimes b \mapsto ab$. Giving structure of $R$-module to $R \otimes_S R$ by multiplying on the left, we define the $R$-module
$$\Omega^{(p)}_{R/S} := I/I^{p+1}.$$

Define the map $d_p^R : R \rightarrow \Omega^{(p)}_{R/S}$,\, $a \mapsto (1 \otimes a - a \otimes 1) + I^{p+1}$. This map is a derivation of order $p$, and its image generates $\Omega^{(p)}_{R/S}$ as an $R$-module (see  \cite[Chapter II-1]{Nakai}).

 The $R$-module $\Omega^{(p)}_{R/S}$ is called the {\it module of K\"ahler differentials of order} $p$ of $R$ over $S$. The map $d_p^R$ is called the canonical derivation of $R$ in $\Omega^{(p)}_{R/S}$. It comes equipped with a universal derivation $d_{S/R}\in {\rm Der}_R^p(S, \Omega^{(p)}_{S/R})$ with the property that composition with $d_{S/R}$ yields an isomorphism $\Hom_R\left(\Omega^{(p)}_{R/S},R\right) \cong \Der_S^p(R)$ (\cite[Proposition 1.6]{O}). To see more properties regarding modules of derivations and K\"ahler differentials, we recommend \cite{O, NKI, Nakai}.

\begin{convention}
It's important to note that the differential module and K\"ahler differentials may not be finitely generated. However, they are finitely generated in certain cases. For instance: If $R$ is essentially of finite type over $S$.
If $S=k$ is a field with a valuation and $R$ is an analytic $k$-algebra, meaning $R$ is finitely generated over a convergent power series ring $k\{x_1,\ldots, x_n\}$.
 If $S=k$ is a field, $(R, \mathfrak{m})$ is a complete local ring, and $R/\mathfrak{m}$ is a finite extension of $k$. For this reason, in this paper, for each $p$, we will consider $\Omega^{(p)}_{R/S}$ to be finitely generated as an $R$-module. In particular, by the universal property,  \({\rm Der}_S^p(R)\) will also be finitely generated as an $R$-module. Thus, in this section, the rings considered are Noetherian local with residue field $k$, and \({\rm Der}_S^p(R)\) will also always be considered a non-zero finitely generated \(R\)-module for each integer \(p \geq 1\). 
\end{convention}


\begin{conjecture}[Zariski-Lipman] \label{conjecture2} Let $R$ be a Noetherian local ring. Then
$R$ is a regular local ring if and only if $\Der_k(R)$ is a free $R$-module.
\end{conjecture}

Following the same direction as the Zariski-Lipman Conjecture, related to the finiteness of the projective dimension of the derivation modules, another important open problem is central in this subject. 

\begin{conjecture}[Herzog-Vasconcelos]
Let $R$ be a Noetherian local ring.  Then $\Der_k(R)$ is a free $R$-module provided ${\rm pd}_R{\rm Der}_k(R) < \infty$.
\end{conjecture}

\begin{conjecture}[Strong Zariski-Lipman Conjecture]\label{conjecture4}
Let $R$ be a Noetherian local ring. Then $R$ is regular provided ${\rm pd}_R{\rm Der}_k(R) < \infty$.
\end{conjecture}
Differently to the Zariski conjecture, the Strong Zariski-Lipman conjecture seems to be widely open, with some exceptions in specific cases (see \cite[Section 4]{He-M} and \cite{HV}).
Due to this open problem, the following  question  introduced by Ludington \cite{Lu} arises.

\begin{questions} (Generalized  Strong  Zariski-Lipman Conjecture) Assume that for some integer $p\geq 1$, \( \operatorname{pd}_{R} \operatorname{Der}_S^p(R) < \infty \). Under what assumptions on \( R \) and $S$, and for which values of \( p \) does this imply that \( R \) is regular?
\end{questions}

The development of the questions in the case \( p > 1 \) has also been conducted, but to a lesser extent compared to the case \( p = 1 \). For further investigations, for example, we refer to Graf \cite{Graf} and Miller-Vassiliadou \cite[Section 4]{Milher}, who study the regularity of an algebraic variety or analytic variety, assuming in the case that modules of differentials or K\"ahler differentials are locally free. It is also worth noting that the conjectures above are generally false for positive characteristic. For example, for \( n = 1 \), it is sufficient to consider \( R = \mathbb{Z}_2[[x,y]]/(x^2-y^2) \) and \( S = \mathbb{Z}_2 \). Additionally, there is also a counterexample for \( n > 1 \), for which we refer to \cite{Lu} for details.
Due to this, from now, our  rings will be equicharacteristic and of characteristic zero. The next result gives a numerical criterion for the true of Generalized  Strong  Zariski-Lipman Conjecture.
 
\begin{theorem}\label{mainteo} Suppose that for some integer $p\geq 1$ and  $n>0$,  $\Der_S^p(R)$ satisfies the inequality
\( \beta_{n}^{R\ltimes k}\left(\Der_S^p(R)\right) \leq \beta_{n-1}^{R\ltimes k}\left(\Der_S^p(R)\right) \). If $\operatorname{Der}_S^p(R)$ is free, then $R$ is regular.
\end{theorem}
\begin{proof} By the proof of Proposition \ref{bettinumbers=}, one obtains that 
 $$\beta_{n}^{R\ltimes k}(\Der_S^p(R))=\sum_{i=0}^{n}\beta_{i}^{R}(\Der_S^p(R))B_{n-i},$$ where $B_{i}=\sum_{j=1}^{i}\beta_{j-1}^{R}(k)B_{i-j}$. As $\Der_S^p(R)$ is free  $R$-module by assumption, we have that $\beta_{i}^{R}(\Der_S^p(R))=0$ for all $i>0$. Then
$$\beta_{n}^{R\ltimes k}(\Der_S^p(R))=\beta_{0}^{R}(\Der_S^p(R))B_{n}=\beta_{0}^{R}(\Der_S^p(R))\left [ B_{n-1}+\sum_{j=2}^{n}\beta_{j-1}^{R}(k)B_{n-j} \right ]$$
and
$$\beta_{n-1}^{R\ltimes k}\left(\Der_S^p(R)\right)=\beta_{0}^{R}\left(\Der_S^p(R)\right)B_{n-1}.$$
By hypothesis, $\beta_{n}^{R\ltimes k}(\Der_S^p(R))\leq\beta_{n-1}^{R\ltimes k}(\Der_S^p(R))$ gives 
$$\sum_{j=2}^{n}\beta_{j-1}^{R}(k)B_{n-j}=0.$$
Since $B_{i}\geq 0$ for all $i\geq 0$  and $B_{0}=1$, one has that $\beta_{n-1}^{R}(k)=0$, that is, the residue field $k$ has finite projective dimension over $R$, and this provides that $R$ is regular.
\end{proof}
As an immediate consequence we derive the following result, specifically directed towards Conjecture \ref{conjecture2}.
\begin{corollary}\label{ZLCclassica}
 Zariski-Lipman's conjecture holds true for any ring $R$ provided for some $n>0$ 
$$\beta_{n}^{R\ltimes k}(\Der_k(R))\leq\beta_{n-1}^{R\ltimes k}(\Der_k(R)).$$
\end{corollary}





\begin{corollary}\label{HVconj} Herzog-Vasconcelos's conjecture implies the Strong Zariski-Lipman Conjecture for any ring $R$ provided for some $n>0$ 
$$\beta_{n}^{R\ltimes k}(\Der_k(R))\leq\beta_{n-1}^{R\ltimes k}(\Der_k(R)).$$
\end{corollary}

\begin{remark}
An alternative way to give a complete answer for the Herzog-Vasconcelos conjecture over any ring $R$, is to show that  $R\ltimes \Der_k(R)$ is a Cohen-Macaulay ring. In fact,  if $R\ltimes \Der_k^{p}(R)$ is Cohen-Macaulay, then $\depth_{R}(\Der_k^{p}(R))\geq \dim (R)$ by Remark \ref{setupremark} (iv). Therefore $\Der_k^{p}(R)$ is maximal Cohen-Macaulay $R$-module. The freeness of  $\Der_k^{p}(R)$  follows by the Auslander-Buchsbaum formula \cite[A.5. Theorem, p. 310]{LW}. This interesting relation changes  completely a homological  problem (to show that  $\Der_k(R)$ is free), to a  structural investigation of a certain ring, i.e., when $R\ltimes \Der_k(R)$ is Cohen-Macaulay.
This allow us to pose  the following question.  

\end{remark}

\begin{question}\label{conjecturanossa}
Let $R$ be a Cohen-Macaulay local ring. If ${\rm pd}_R{\rm Der}_k(R) < \infty$, is it true that  $R\ltimes \Der_k(R)$ must be   a  Cohen-Macaulay ring? 

\end{question}

\begin{corollary}\label{HVconj2} If Question \ref{conjecturanossa} is true, then Herzog-Vasconcelos's conjecture is true. 
\end{corollary}

\begin{example}
The following example illustrates that, in positive characteristic, the Cohen-Macaulay hypothesis on the ring $R$ can not be dropped in Question \ref{conjecturanossa}.
Let $R := \mathbb{F}_2[[X,Y]]/(X^4,X^2Y^2)$. For simplicity, we set $k := \mathbb{F}_2$, $f_1 := X^4$ and $f_2 := X^2Y^2$. For any $f \in k[[X, Y]]$, consider $df := \frac{\partial f}{\partial X}dx + \frac{\partial f}{\partial Y}dy$. Note that $\Omega_{R/k}^{(1)} \cong Rdx \oplus Rdy / (df_1, df_2)$, where $Rdx \oplus Rdy$ is a free $R$-module with base $\{dx, dy\}$. Also, one has $df_1=4X^3dx=0$ and $df_2=2Y^2Xdx + 2YX^2dy=0.$ Thus, $\Omega_{R/k}^{(1)}$ is a free $R$-module of rank 2. Note that $\text{Der}_k(R) \cong \text{Hom}_R\left(\Omega_{R/k}^{(1)}, R\right) \cong R^2$. This yields
$
\text{pd}_R(\text{Der}_k(R)) = 0 < \infty.
$
Now,  Remark $\ref{setupremark}(ii)-(iii)$ provides that $\dim(R\ltimes \text{Der}_k(R))=\dim(R)$ and $\depth{(R\ltimes \text{Der}_k(R))}=\depth{R}$. But, as $R$ is not Cohen-Macaulay, one obtains that $R \ltimes \text{Der}_k(R)$ is also not Cohen-Macaulay ring.
\end{example}

\bigskip



\noindent{\bf Acknowledgements.} 
During the development of this work, we received the sad news that Professor J\"urgen Herzog passed away. The authors had the honor of meeting Professor Herzog during the International Meeting in Commutative Algebra and its Related Areas (IMCARA-2017) - S\~ao Carlos - Brazil. His generosity and valuable mathematical contributions will certainly last through generations. The authors would like to thank Saeed Nasseh and Souvik Dey for their valuable suggestions for this paper.

\end{document}